\documentclass[11pt,reqno]{amsart}
\usepackage[]{amsmath,amssymb,amsfonts,latexsym,amsthm,enumerate}
\usepackage{amssymb,bbm,mathrsfs,tikz}
\usepackage[normalem]{ulem}
\usepackage{enumerate,graphicx,paralist}
\usepackage[font=small]{caption}
\usepackage[noadjust]{cite}

\numberwithin{equation}{section}
\usepackage[colorlinks=true, pdfstartview=FitV, linkcolor=blue,
  citecolor=blue, urlcolor=blue,pagebackref=false]{hyperref}
\newcommand{\red}[1]{
#1
}
\newcommand{\green}[1]{{\color{green} #1}}

\usepackage[colorinlistoftodos]{todonotes}
\presetkeys{todonotes}{inline, color=green}{}
\usepackage{comment}

\newtheorem{theorem}{Theorem}[section]
\newtheorem*{theorem*}{Theorem}
\newtheorem{lemma}[theorem]{Lemma}

\newtheorem{remark}[theorem]{Remark}
\newtheorem*{remark*}{Remark}

\theoremstyle{definition}{

\newtheorem*{definition*}{Definition}

\newtheorem*{question*}{Question}
\newtheorem*{example*}{Example}
\newtheorem*{examples*}{Examples}
}

\newcommand{\abbr}[1]{{\sc{\lowercase{#1}}}}

\newcommand{\A}{\mathbb A}

\newcommand{\E}{\mathbb E}
\newcommand{\N}{\mathbb N}
\renewcommand{\P}{\mathbb P}

\newcommand{\R}{\mathbb R}

\newcommand{\cK}{\mathcal K}
\newcommand{\cA}{{\mathcal A}}
\newcommand{\cB}{{\mathcal B}}

\newcommand{\cF}{{\mathcal F}}

\newcommand{\cZ}{{\mathcal Z}}

\newcommand{\An}[1][n]{\A_{#1}^+}
\newcommand{\bAn}[1][n]{\bar\A_{#1}^+}
\newcommand{\Ltt}{L^2(\red{\An \times \An})}

\newcommand{\bE}{\mathbf E}
\newcommand{\uh}{\underline h}
\newcommand{\un}{\underline n}
\newcommand{\uu}{\underline u}
\newcommand{\uv}{\underline v}
\newcommand{\uw}{\underline w}
\newcommand{\ux}{\underline x}
\newcommand{\uy}{\underline y}
\newcommand{\uz}{\underline z}
\renewcommand{\d}{\mathrm{d}}
\newcommand{\bcfree}{\mathfrak{f}}
\newcommand{\bczero}{\mathfrak{0}}
\newcommand{\ct}{t}
\newcommand{\cs}{s}
\newcommand{\crr}{r}

\newcommand{\fa}{\mathfrak a}
\newcommand{\fb}{\mathfrak b}
\newcommand{\fc}{\mathfrak c}

\newcommand{\fB}{\mathfrak B}

\newcommand{\sC}{\mathsf C}

\newcommand{\Dim}{\textsc{d} }

\newcommand{\one}{\mathbbm{1}}
\renewcommand{\epsilon}{\varepsilon}

\author{Amir Dembo}
\address{Amir Dembo\hfill\break
Mathematics  Department and Statistics Department\\ Stanford University\\
Stanford, CA 94305, USA.}
\email{adembo@stanford.edu}

\author{Eyal Lubetzky}
\address{Eyal Lubetzky\hfill\break
Courant Institute %of Mathematical Sciences
\\ New York University\\
251 Mercer Street\\ New York, NY 10012, USA.}
\email{eyal@courant.nyu.edu}

\author{Ofer Zeitouni}
\address{Ofer Zeitouni\hfill\break
Department of Mathematics\\
Weizmann Institute of Science\\
Rehovot 76100, Israel\\
and
Courant Institute\\
New York University\\
251 Mercer Street\\ New York, NY 10012, USA.}
\email{ofer.zeitouni@weizmann.ac.il}

\title[Limit of line ensembles with geometric area tilt]{On the limiting law of line ensembles of Brownian~polymers with geometric area tilts}

\dedicatory{Dedicated to the memory of Dima Ioffe}

\begin{document}

\begin{abstract}
We study the line ensembles of non-crossing Brownian bridges above a hard wall, each tilted by the area of the region below it with geometrically growing pre-factors. This model, which mimics the level lines of the $(2+1)$\Dim \abbr{SOS} model above a hard wall, was studied in two works from 2019 by Caputo, Ioffe and Wachtel. In those works, the tightness of the law of the top $k$ paths, for any fixed $k$, was established under either zero or free boundary conditions, which in the former setting implied the existence of a limit via a monotonicity argument. Here we address the open problem of existence of  a limit under free boundary conditions: we prove that as the interval length, followed by the number of paths, go to $\infty$, the top $k$ paths converge to the same limit as in the zero boundary case, as conjectured by Caputo, Ioffe and Wachtel.\\[0.6em]
\iffalse
\noindent
\textsc{R\'{e}sum\'{e}.} Nous \'{e}tudions l'ensemble de lignes d\'{e}termin\'{e} par des mouvements Browniens non-intersectant au-dessus d'un mur solide. Ce mod\`{e}le, qui imite les lignes de niveaux
du mod\`{e}le  $(2+1)$\Dim \abbr{SOS} au-dessus d'un mur, a \'{e}t\'{e} \'{e}tudi\'{e} en 2019 par  Caputo, Ioffe et Wachtel. Dans ces travaux, la tension de la loi des $k$ lignes hautes, pour chaque $k$ fixe, a \'{e}t\'{e} obtenue sous des condition nulles au bord ou des conditions libres au bord. Dans le premier cas, ca implique l'existence d'une limite par un argument de monotonicit\'{e}.
Nous adressons ici le probl\`{e}me ouvert d'existence d'une limite sous des conditions libres au bord: nous d\'{e}montrons que quand la longueur de l'intervalle, suivi par le nombre de lignes, tend vers l'infinie, les $k$ lignes hautes convergent vers la m\^{e}me limite que dans le  cas de conditions nulles au bord, comme conjectur\'{e} par
 Caputo, Ioffe et Wachtel.
\fi
\end{abstract}

%{\mbox{}\vspace{-0.6cm}
\maketitle
%}
%\vspace{-0.6cm}

\section{Introduction}\label{sec:intro}
Entropic repulsion in low temperature $(2+1)$\Dim crystals above a hard wall has been the subject of extensive study in statistical physics. Whereas in the absence of a wall, the surface of the crystal would typically be rigid at height $0$, in the presence of a wall, the surface is propelled in order to increase its entropy (i.e., to allow thermal fluctuations going downward), and becomes rigid at some height level which diverges with the side length $L$ of the box.

A rigorous study of this phenomenon in the $(2+1)$\Dim Solid-On-Solid \abbr{(SOS)} model---a low temperature approximation of the 3\Dim Ising model---dates back to Bricmont, El Mellouki and Fr\"{o}hlich~\cite{BEF86} in 1986, where it was shown that, in the presence of a hard wall at height $0$, the typical height of a site in the bulk is propelled to order $\log L$. Thereafter, a detailed description of the shape of this random surface was obtained by Caputo et al.~\cite{CLMST12,CLMST14,CLMST16}, showing that it typically becomes rigid at a height which is one of two consecutive (explicit) integers, through a sequence of nested level lines each encompassing a $(1-\epsilon)$-fraction of the sites (analogous behavior was later established~\cite{LMS16} for the more general family of $|\nabla\phi|^p$-random surface model, where the \abbr{SOS} model is the case $p=1$).
The level lines near the center sides of the box behave as random walks---a ubiquitous feature of interfaces in low temperature spin systems---albeit with cube-root fluctuations, as their laws are tilted by the entropic repulsion effect. The lower the level line, the higher the reward is for generating spikes going downward, and as such, the tilting effect of the level lines increases exponentially as the height decreases.

Whereas the $2$\Dim Ising model with a pinning potential is known~\cite{IOSV21} to have an interface converging to a Ferrari--Spohn diffusion, the behavior in the \abbr{SOS} model---where there are $H\asymp\log L$ interacting level lines, each constrained not to cross its neighbors and inducing a tilt which is a function of the area it encompasses and its height---is far from being understood (see the review~\cite{IV18} for more information).

In this work, we investigate the limiting law of a line ensemble that was
 studied by Caputo, Ioffe and Wachtel~\cite{CIW19a,CIW19b} to model the level lines of the \abbr{SOS} model in the presence of a hard wall: each level line, $X_1,X_2,\ldots$, where $X_1$ is the top one, is tilted by the area below it, with the coefficients of these area tilts increasing geometrically.

For more perspective on this model in the context of other models of Brownian polymers constrained above a barrier, starting from the influential work of Ferrari and Spohn~\cite{FS05} (the model there being equivalent by Girsanov's transformation---cf.~\cite{MZ16}---to a Brownian excursion with an area tilt), see, e.g.,~\cite{CIW19a,ISV15,IVW18} and the references therein.

Define
\[ \An = \{ \ux\in \R^n \,:\; x_1 > x_2 > \ldots > x_n > 0\}\,,\]
its closure $\bAn$ and, for a designated interval
\[I=[\ell,r] \qquad (\ell<r\in \R)\,,\]
let
\[ \Omega_{n}^{I} = \left\{ X\in \sC(I; \mathbb{R}^n)
\,:\; X(t)\in \An\mbox{ for all $t\in I$}\right\}\,.\]
(Here, for $T\subset \mathbb{R}$ and $\mathcal{X}$ a topological space, we denote by $\sC(T;\mathcal{X})$ the space of continuous functions from $T$ to $\mathcal{X}$, equipped with the topology of uniform convergence on compact subsets of $\mathbb{R}$).
Further define the area tilt of $Y\in \sC(I;\mathbb{R})$ to be
\[ \cA_I(Y) = \int_I Y(t)\, \d t\,,\]
and, for given tilt parameters $\fa>0$ and $\fb>1$ and
endpoints $\ux=(x_1,\ldots,x_n)\in \An$ and $\uy=(y_1,\ldots,y_n)\in \An$, the partition function
\begin{equation}\label{eq:partition-func}
Z_{n}^{\ux,\uy,I} = \bE_{n}^{\ux, \uy,I} \left[\one_{\Omega_{n}^{I}}\, e^{-\fa \sum_{i=1}^n  \fb^{i-1} \cA_I(X_i(\cdot))}\right]\,,\end{equation}
in which $\bE_n^{\ux,\uy,I} = \bigotimes_{i=1}^n \bE_1^{x_i,y_i,I}$ and the expectation $\bE_1^{x,y,I}$ for $I=[\ell,r]$ is w.r.t.\ the (unnormalized) path measures of the Brownian bridge which starts at $x$ at time $\ell$ and ends at $y$ at time~$r$;
that is, the total mass of $\bE_1^{x_i,y_i,I}$ is $\phi_{r-\ell}(y_i-x_i)$,
where hereafter $\|\cdot\|$ stands for the relevant Euclidean norm, with
\begin{equation}\label{def:Gauss}
\phi_v(\ux):=(2\pi v)^{-k/2} e^{-\|\ux\|^2/(2v)}
\end{equation}
denoting the density of a centered Gaussian vector
of independent coordinates of variance $v$, whose dimension $k$ is
implicitly given by the argument we use.
% $(1/\sqrt{2\pi(r-\ell)}) \exp(-\frac{(y_i-x_i)^2}{2(r-\ell)})$.
(At no point in our analysis will we need to adjust the tilt parameters $(\fa,\fb)$, and as such we do not include them in the notation for brevity.)

Let $\cB_n=\cB_{n,I}$ be the Borel $\sigma$-field on $\sC(I,\R^n)$. (We omit $I$ from the notation when no confusion occurs.)  For $\Gamma\in\cB_{n,I}$ define
\begin{equation}\label{eq:P-ux-uy-def} \P_{n}^{\ux,\uy,I}(\Gamma) := \frac1{Z_{n}^{\ux,\uy,I}} \bE_n^{\ux,\uy,I}\left[\one_\Gamma \one_{\Omega_{n}^I}\, e^{-\fa \sum_{i=1}^n \fb^{i-1} \cA_I(X_i(\cdot))}\right]\,.
\end{equation}
Consider $I_T=[-T,T]$. Two classes of boundary conditions that are of interest are:
\begin{enumerate}[(a)]
\item \emph{Zero boundary conditions}: fixing both $\ux$ and $\uy$ to be zero:
\[ \mu^\bczero_{n,T} = \P_{n}^{\underline 0,\underline 0, I_T} \;\]
(more precisely, this is the limit of
$\P_{n}^{\epsilon \ux,\epsilon \uy, I_T}$ as $\epsilon \downarrow 0$,
which by stochastic domination exists and is independent of the
fixed $\ux$, $\uy$ in $\An$ which one uses).
\item
\emph{Free boundary conditions} with respect to a $\sigma$-finite measure:
averaging $\bE_n^{\ux,\uy,I_T}[\cdot]$
over $\ux,\uy$ according to a specified
$\sigma$-finite measure $\Theta_n$ on~$\R^n$:
% \blue{(e.g. Lebesgue):}
\[ \mu^\bcfree_{n,T}(\Gamma) =
\frac1{\cZ_{n,T}^{\bcfree}} \int_{\An} \int_{\An} \bE_n^{\ux,\uy,I_T}\left[\one_\Gamma \one_{\Omega_n^{I_T}}\, e^{-\fa \sum_{i=1}^n \fb^{i-1} \cA_{I_T}(X_i(\cdot))}\right] \,\Theta_n(\d \ux) \Theta_n(\d \uy) \,,
\]
where
\[ \cZ_{n,T}^\bcfree := \int_{\An} \int_{\An}
\bE_n^{\ux,\uy,I_T}\left[ \one_{\Omega_n^{I_T}}\, e^{-\fa\sum_{i=1}^n \fb^{i-1} \cA_{I_T}(X_i(\cdot))}\right] \,
\Theta_n(\d \ux) \Theta_n(\d \uy) \,.
\]
\end{enumerate}
We refer to such $\mu^{\bcfree}_{n,T}$ as $\Theta_n$-free boundary conditions, reserving
Leb-free for the special case of Lebesgue $\Theta_n$, considered in \cite{CIW19a,CIW19b}.
Caputo, Ioffe and Wachtel show in~\cite{CIW19a,CIW19b}
that $\mu_{n,T}^\bczero$ converges to a limit $\mu_\infty^\bczero$ as $n,T\to\infty$
(moreover, they proved that for any fixed $n$, the measures $\mu_{n,T}^\bczero$ converge as $T\to\infty$ to a limit $\mu_n^\bczero$, which then converges to $\mu_\infty^\bczero$ as $n\to\infty$),
and that for any $c>0$, the family of Leb-free distributions
$\{\mu_{n,T}^\bcfree\}_{n\geq 1, T>c}$ is tight.
In this and subsequent statements, the measure $\mu_\infty^\bczero$ is defined on
$\sC(\mathbb{R},\mathbb{R}^{\N})$, and the convergence is in the sense that for any compact set $\cK\subset \mathbb{R}$, integer $k\in \N$ and fixed bounded, continuous
function
$f: \sC(\cK;\mathbb{R}^k) \mapsto \mathbb{R}$, we have that
\begin{equation}
    \label{def:conv}
\lim_{n,T\to\infty} \int f(X_1,\ldots,X_k) \mu_{n,T}^\bczero(\d X)= \int f(X_1,\ldots,X_k)\mu_\infty^\bczero (\d X)\,.
\end{equation}
For Leb-free boundary conditions, Caputo, Ioffe and Wachtel
conjectured that $\mu_{n,T}^\bcfree$ converges as well, and to the same limit
$\mu^\bczero_\infty$ as $n,T\to\infty$.
Our main result confirms that when $T\to\infty$ followed by~$n\to\infty$, this holds
more generally, whenever for $n \ge 1$,
\begin{align}
\fc_n  &:= \limsup_{r \to \infty} r^{-1} \log \Theta_n(\An \cap \{x_1 \le r\}) < \infty
\label{eq:exp-growth}
\end{align}
(in particular, note that $\fc_n = 0$ when $\Theta_n$ is Lebesgue measure).

\begin{theorem}\label{thm:1}
Assuming \eqref{eq:exp-growth}, for
any fixed tilt parameters $\fa>0$ and $\fb>1$ and any fixed integer $n$, the measures $\mu_{n,T}^\bcfree$ and $\mu_{n,T}^\bczero$ have the same weak limit as $T\to\infty$. In particular, if we denote by $\mu_\infty^\bczero$ the limit of $\mu_{n,T}^\bczero$ as $n,T\to\infty$, then
\[ \exists\lim_{n\to\infty} \lim_{T\to\infty} \mu_{n,T}^\bcfree =\mu_\infty^\bczero\,.
\]
\end{theorem}

\begin{remark} Our proof easily extends to allow in
$\mu_{n,T}^\bcfree$ different measures for $\ux$ and for $\uy$ (as long as
both satisfy \eqref{eq:exp-growth}). Note that Theorem \ref{thm:1} is optimal in terms of $\Theta_n$,
as merely having $\cZ_{n,T}^\bcfree$ finite, requires that $\fa \, T \ge \fc_n$ (see \eqref{eq:mu-free-1}
and \eqref{eq:psi-bd} at $s=T$). Further, $\sup_n \{\fc_n\}$ must be finite if aiming to exchange
the order of limits in $n$ and $T$.
\end{remark}

Our proof of Theorem~\ref{thm:1} employs the Markovian structure of the problem. In a first step we introduce a (sub-)Markovian  Kernel $K_t$, see~\eqref{eq-kernel}. The key part of the proof is Lemma~\ref{lem:K1-compact}, where we prove that $K_1$ is compact in the appropriate $L^2$ space; the proof of the lemma involves probabilistic arguments. With the lemma, standard contraction arguments, detailed in Section~\ref{subsec-2.1}, yield the exponential decay (in $T$) of the dependence in the boundary conditions. We note that some care is needed here due to the non-compactness of the set of possible boundary conditions, but that non-compactness was already handled in~\cite{CIW19a}.

Many interesting open questions remain, chief among which, perhaps, is describing the
limiting process $X_\infty(\cdot)$ (on, say, the interval $[0,1]$). We refer to~\cite{CIW19a,CIW19b} for a list of such problems 
and note in passing that from \eqref{eq:mu-free-4} and
\abbr{PDE} theory, one can verify that
% one finds that
the
%weak
limit of  $\mu^\bcfree_{n,T}$ when $T \to \infty$,
is the stationary solution of the
Langevin \abbr{SDE} for invariant probability density $\varphi_1^2$ on $\An$ (where
the Perron-Frobenius eigenvector $\varphi_1$ of $K_1$ is the positive $C^2$-solution of
the elliptic \abbr{PDE}
%$u(s,\cdot) = K_s \varphi_1 (\cdot)$ is the unique
%$C^{1,2}$-solution of the
% heat equation parabolic \abbr{PDE}
%$u_t=$\frac{1}{2} \Delta u - \fa \langle \underline{\fb},\uz\rangle u$
$\frac{1}{2} \Delta u = (c + \fa \langle \underline{\fb},\ux\rangle) u$
with
%a linear in $\ux \in \An$ killing rate
%on $(0,\infty) \times \An$ with
%and
Dirichlet boundary
conditions at $\partial \An$ and the largest possible $c
%=\log \lambda_1
<0$).
% $\bAn \setminus \An$ starting
%and positive initial condition
%$
% u(0,\cdot)= \varphi_1$. In particular, such $u$ is strictly positive on $\R_+ \times \An$. Further,
%with $\fc_n$ of \eqref{eq:exp-growth} finite, fixing $r$ finite such that
%$\Theta_n(\An \cap \{x_1 \le r\})>0$, for any $s \ge \ell$,

\section{Proof of main result}
Fix the tilt parameters $\fa>0$ and $\fb>1$, and let $n\geq 1$ be an integer.
Throughout this proof, for $X\in \Omega_n^I$,
we use the abbreviated notation
\[ \cA_I(X(\cdot)) := \fa \sum_{i=1}^n \fb^{i-1} \cA_I(X_i(\cdot))\,.\]
Let $\Gamma\in\cB_{n,[0,1]}$, and define
\[ K^\Gamma_1(\ux,\uy) = \bE_n^{\ux,\uy,[0,1]}\left[\one_\Gamma \one_{\Omega_n^{[0,1]}}\, e^{-\cA_{[0,1]}(X(\cdot))}\right]\,,
\]
which we view as a linear operator on $L^2(\An)=L^2(\An,{\rm Leb})$:
\[ (K^\Gamma_1 f)(\ux) = \int_{\An} K^\Gamma_1(\ux,\uy) f(\uy) \,\d \uy\,.\]
With a slight abuse of notation,
we continue to write $K_1^\Gamma$ also when $\Gamma\in \cB_{n,\mathbb{R}}$, in which case we understand that $\Gamma$ was replaced by its restriction to the interval $[0,1]$. With this convention in mind,
we will further be interested in the semigroup
\begin{equation}
    \label{eq-kernel}
K^\Gamma_t(\ux,\uy) = \bE_n^{\ux,\uy,[0,t]}\left[\one_\Gamma \one_{\Omega_n^{[0,t]}}\, e^{-\cA_{[0,t]}(X(\cdot))}\right]\,.
\end{equation}
When referring to the case $\Gamma=\Omega_n^\R$ (i.e., the indicator $\one_\Gamma$ within the expectation in the definition of $K_1^\Gamma$ is omitted), we simply write $K_t$ (with no superscript) in lieu of~$K_t^{\Omega_n^\R}$, noting that $K_t(\ux,\uy)$ is precisely the partition function $Z_n^{\ux,\uy,[0,t]}$ from~\eqref{eq:partition-func}.

Observe that $K_t$ is symmetric, in that $K_t(\ux,\uy)=K_t(\uy,\ux)$,
as well as positivity preserving:
\[ (K_t f)(\ux) =\int K_t(\ux,\uy)f(\uy)\d\uy \geq 0\qquad\mbox{whenever}\qquad f \geq 0\,.\]
As $K_t$ is symmetric, and given by a continuous time Markov process with killing, it is positive definite (this follows, e.g., by~\cite[Theorems~1.3.1, Lemma~1.3.2 and Theorem~6.1.1]{Fukushima}, all applied to Example 1.2.3 there).
A key ingredient in the proof will be that $K_1$ is furthermore relatively compact:
\begin{lemma}\label{lem:K1-compact}
For every fixed $n$, the range of the symmetric positive definite operator
%$K_1$ given by
\[ (K_1 f)(\ux) = \int_{\An} \bE_n^{\ux,\uy,[0,1]}\left[ \one_{\Omega_n^{[0,1]}}\, e^{-\cA_{[0,1]}(X(\cdot))}\right]f(\uy)\,\d\uy
\]
consists of continuous functions on $\An$, and the operator $K_1$ is compact w.r.t.\ $L^2(\An)$.
\end{lemma}

\subsection{Proof of Theorem~\ref{thm:1} modulo Lemma~\ref{lem:K1-compact}}
\label{subsec-2.1}
We consider throughout the convergence over the interval $[0,1]$, the changes
needed for considering other compact sets (as the set $\cK$ in \eqref{def:conv}) are minimal. Aiming
to express the measures $\mu_{n,T}^\bczero$ and $\mu_{n,T}^\bcfree$ in terms of the operator $K_t$, we
define for $s >0$
\[
\psi_{\cs}(\uu) = \int K_{\cs}(\uu,\ux)\, \Theta_n(\d \ux) \,.
\]
%noting that $\psi_{s}(\uu)$ is well-defined (albeit a-priori possibly infinite) as $K_s$ is non-negative.
 Setting $s>0$ large enough so that $\psi_s \in L^2(\An)$,
%\blue{$\fa s > \fc_n$ of \eqref{eq:exp-growth}},
in view of
the symmetry of $K_{s}$ and the semigroup property, our
goal is then to show that for every $\Gamma\in\cB_{n,[0,1]}$, the limit of
\begin{equation}
    \label{eq:mu-free-1}
    \mu_{n,T}^\bcfree(\Gamma) = \frac{\iint \psi_s(\ux) K_{T-\cs}(\ux , \uu)
    K^\Gamma_1(\uu, \uv) K_{T-1-\cs} (\uv, \uy) \psi_s(\uy) \, \d \uu \d \uv \, \d \ux \d \uy} {\iint \psi_s(\ux) \, K_{2T-2\cs}(\ux, \uy)\, \psi_s(\uy) \d\ux \d\uy}
\end{equation}
as $T\to\infty$ exists and coincides with that of
\begin{equation}
    \label{eq:mu-zero-1}
 \mu_{n,T}^\bczero(\Gamma) = \lim_{\epsilon\downarrow 0}\frac{\iint K_{T}(\epsilon \ux, \uu) K_1^\Gamma(\uu, \uv) K_{T-1}(\epsilon \uy, \uv) \,\d \uu \d \uv }{K_{2T}(\epsilon \ux,\epsilon \uy)}\,,
\end{equation}
where by~\cite[Lemma~2.2]{CIW19b}, the limit \eqref{eq:mu-zero-1} exists and
is independent of
%the choice of
$\ux,\uy$ in $\An$.

By the spectral decomposition theorem, the compact positive definite operator~$K_1$ has a discrete spectrum (except for a possible accumulation point at $0$), with positive eigenvalues $\{\lambda_i\}$
and eigenvectors $\{\varphi_i\}$ that form a complete orthonormal basis of $L^2(\An)$
(see, e.g.,~\cite[Thms.~VI.15 and VI.16]{ReedSimonI}). In particular,
\[
    K_1(\ux,\uy) = \sum_{i=1}^\infty \lambda_i \varphi_i(\ux) \varphi_i(\uy) \quad\mbox{for a complete basis $\{\varphi_i\}_{i\geq 1}$ with}\; \left<\varphi_i,\varphi_j\right>_{L^2(\An)}=\delta_{ij}\,.
    %K_1\phi_i &= \lambda_i\phi_i\qquad(i=1,2,\ldots)\,.
\]
With $K_1(\ux,\uy)>0$ throughout $\An \times \An$ (e.g., due to parabolic regularity),
by the generalized Perron-Frobenius Theorem (see, e.g., the version of the
Krein-Rutman Theorem given in [9, Thm. XIII.43]), the
top eigenvalue $\lambda_1$ has a one-dimensional eigen-space and we may choose
the continuous function $\varphi_1$ to be strictly positive on $\An$. That is,
\[ \varphi_1 > 0\qquad\mbox{and}\qquad \lambda_1>\lambda_2\geq \lambda_3 \geq \ldots \geq 0\,.\]
Further, for
any $\crr \ge 1$ and $\ux,\uy\in \An$,
\begin{align}
    \label{eq-K1K2}
    K_{\crr}(\ux,\uy) = \bE_n^{\ux,\uy,[0,\crr]}\left[ \one_{\Omega_n^{[0,\crr]}}\, e^{-\cA_{[0,\crr]}(X(\cdot))}\right]
    &\leq \bE_n^{\ux,\uy,[0,\crr]}\left[ \one_{\Omega_n^{[0,\crr]}}\, e^{-\fa \cA_{[0,\crr]}(X_1(\cdot))}\right] \nonumber \\
    &\le
    \bE_1^{x_1,y_1,[0,\crr]}\left[
    %\one_{\Omega_{1}^{[0,\crr]}}\,
    e^{-\fa\cA_{[0,\crr]}(X_1(\cdot))}\right]\nonumber\\
    & {\;  = \;} \phi_{\crr}(y_1-x_1) e^{-\fa \, \crr \frac{x_1+y_1}2} \E \big[e^{-\fa  \int_0^{\crr} B_s \,\d s}\big]
    \nonumber\\
    &
    \le e^{C_r} 
    %e^{-\frac{(x_1-y_1)^2}{2r}}
    e^{-\fa \, \crr \frac{x_1+y_1}2},
\end{align}
where $\{B_s, s\in [0,\crr]\}$ is
the standard Brownian bridge over $[0,\crr]$ starting and ending at $0$
%, with expectation denoted $\bE$
and $C_r =\frac{\fa^2}{2} \E(\int_0^r B_s \d s)^2$
(using in the second line that the total mass of
$\bE_1^{x_i,y_i,[0,\crr]}$, $i \ge 2$, is at most one, while for the third line recall that
a Brownian bridge between fixed points has the law of the standard bridge
plus a straight line connecting these points). Since $K_{\crr}(\ux,\uy)$ vanishes if either
$\ux\not\in \An$ or $\uy\not\in \An$, it follows that
\begin{align}
    \label{eq-K2bound}
    \iint K_{\crr}(\ux,\uy) \,\d \ux \d\uy &\leq  e^{C_r} \Big[
    %\frac{1}{(n-1)!}
    \int_0^\infty x_1^{n-1}
    e^{-\fa \, r x_1/2} \,\d x_1 \Big]^2 < \infty \,,
    \end{align}
    and
    \begin{align}    \label{eq-K1bound}
    \int K_{\crr} (\ux,\ux) \,\d \ux &\leq e^{C_r} \int_0^\infty x_1^{n-1}   e^{-\fa \, \crr x_1}\,\d x_1
    < \infty \,.
\end{align}
%\blue{are bounded, uniformly over $[1/(2\pi),\tau]$ for any $\tau$ finite.}
Similarly,
by the symmetry of $K_s$, the semigroup property and \eqref{eq-K1K2},
\begin{equation}\label{eq:psi-bd}
\int \psi_{\cs}(\uu)^2\,\d\uu =
%\iint K_s(\ux,\uu) K_s(\uu,\uy) \, \blue{\Theta_n(\d \ux) \Theta_n(\d \uy)} \d\uu =
\iint K_{2\cs}(\ux,\uy) \Theta_n(\d \ux) \Theta_n(\d \uy) \le
e^{C_{2s}} \Big[ \int_{\An} e^{-\fa s x_1} \Theta_n(\d \ux) \Big]^2 <\infty \,,
\end{equation}
provided that $\fa \, s > \fc_n$ of \eqref{eq:exp-growth}, in which case we can
decompose
\begin{equation}\label{eq:l2-dec}
\psi_{\cs} = \sum_{i=1}^\infty \alpha_{i,\cs} \varphi_i\quad\mbox{where}\quad \alpha_{i,\cs} := \left<\psi_{\cs},\varphi_i\right>_{L^2(\An)}\,, \quad \sum_{i=1}^\infty \alpha_{i,\cs}^2  = \|\psi_{\cs}\|_2^2
\end{equation}
(hereafter $\|\cdot\|_2$ denotes the
$L^2(\mathbb{R}^n,{\rm Leb})$-norm, using $\|\cdot\|_{L^2(\An)}$ when restricting the domain
to $\An$). Fixing an integer $\ell > \fc_n/\fa$, we have from \eqref{eq:psi-bd} that
$\{\psi_s\}_{s \in [\ell,\ell+1]}$ is bounded in $L^2(\An)$, and we split any
$T \ge \ell+1$ as $T=t+s$ for $s
%=\ell+\{T\}
\in [\ell,\ell+1)$ and integer $t \ge 1$,
to get from the decomposition~\eqref{eq:l2-dec} that
\begin{align*}
    \int K_{t-1}(\uu,\uy) \psi_{\cs} (\uy) \,\d\uy
   = \sum_{i=1}^\infty \lambda_i^{t-1}\alpha_{i,\cs} \varphi_i(\uu) \,.
\end{align*}
Similarly, for $t\geq 1$,
\begin{align*}
 \iint K_{\ct}(\ux,\uy) \psi_{\cs}(\ux)\psi_{\cs}(\uy)\,\d\ux\d\uy
= \sum_{i=1}^\infty \lambda_i^{\ct}\alpha_{i,\cs}^2 := c_{t,s}
\,.\end{align*}
Hence,~\eqref{eq:mu-free-1} translates for $s=\ell+\{T\}$ and $t=T-s$, into
\begin{align} \mu_{n,T}^\bcfree(\Gamma) &= \frac{1}{c_{2t,s}} \iint \sum_{i=1}^\infty\sum_{j=1}^\infty
\alpha_{i,\cs} \alpha_{j,\cs} \lambda_i^{\ct} \lambda_j^{\ct-1} \varphi_i(\uu) \, K^\Gamma_1(\uu,\uv)\, \varphi_j(\uv) \,\d\uu\d\uv \,.
\label{eq:mu-free-2}
\end{align}
%Note that the numerator features the indicator $\one_\Gamma$ through the middle term $K_1^\Gamma(\uu,\uv)$.
Looking at $K^\Gamma_1$ and arguing as we did for $K_1$, we see that for any $\uu\in\R^n$,
\[ \int K^\Gamma_1(\uu,\uv)^2\, \d\uv \leq \int K_1(\uu,\uv)^2\,\d\uv = K_2(\uu,\uu) < \infty \,,\]
where the equality holds by the symmetry of $K_1$ and the definition of $K_t$ and the last inequality by~\eqref{eq-K1K2}.
In other words, $K_1^\Gamma(\uu,\cdot)\in L^2(\An)$ for every $\uu\in\R^n$.
Moreover, by
\eqref{eq-K1bound}
we have that
\[ \iint K^\Gamma_1(\uu,\uv)^2 \,\d\uu\d\uv \leq \iint K_1(\uu,\uv)^2\,\d\uu\d\uv =\int K_2(\uu,\uu) \d\uu<\infty\]
and it follows that
\[ K_1^\Gamma \in \Ltt \,.\]
A complete orthonormal system $\{\varphi_i\}$ w.r.t.\ $L^2(\An)$ induces a complete orthonormal system $\{\varphi_i\otimes\varphi_j\}_{i,j\geq 1}$ w.r.t.\ $\Ltt$; hence, we may decompose $K^\Gamma_1$ into
\[ K^\Gamma_1(\uu,\uv) = \sum_{i,j}\gamma_{i,j}\varphi_i(\uu)\varphi_j(\uv)\]
where
\[ \gamma_{i,j} := \iint  K^\Gamma_1(\uu,\uv) \varphi_i(\uu) \varphi_j(\uv)\,
\d\uu\d\uv \,, \qquad
\sum_{i,j \geq 1} \gamma_{i,j}^2 =\| K^\Gamma_1 \|^2_{\Ltt}
< \infty\,.
\]
This reduces~\eqref{eq:mu-free-2} into $\mu_{n,T}^\bcfree(\Gamma) = \Xi^{(1)}_{n,T}/\Xi^{(2)}_{n,T}$ where
\begin{equation}
\label{eq:mu-free-3}
%\mu_{n,T}^\bcfree(\Gamma) = \frac{\Xi^{(1)}_{n,T}}{\Xi^{(2)}_{n,T}} \,, \quad
\Xi^{(1)}_{n,T} := \sum_{i,j \ge 1} \gamma_{i,j} \alpha_{i,\cs} \alpha_{j,\cs} \widehat\lambda_i^{\ct} \widehat\lambda_j^{\ct-1} \,,
\qquad
\Xi^{(2)}_{n,T} :=
\lambda_1\sum_{i=1}^\infty \widehat\lambda_i^{2\ct}\alpha_{i,\cs}^2 \,,
\end{equation}
and the rescaled eigenvalues $\widehat\lambda_i := \lambda_i / \lambda_1 \in [0,1]$ $(i=1,2,\ldots)$ satisfy
\[ \widehat\lambda_i=1\qquad\mbox{and}\qquad\sup_{i>1}\widehat\lambda_i \leq 1-\delta\quad\mbox{for $\delta=(\lambda_1-\lambda_2)/\lambda_1>0$}\,.
\]
We immediately see that $\Xi_{n,T}^{(2)}$ of ~\eqref{eq:mu-free-3} satisfies
\begin{equation}\label{eq:err-Xi2}
\lambda_1\alpha_{1,\cs}^2  \leq \Xi_{n,T}^{(2)} \leq
\lambda_1\alpha_{1,\cs}^2 + \lambda_1 (1-\delta)^{2\ct}\|\psi_{\cs}\|_2^2\,.
\end{equation}
Further, with $K_s \varphi_1 = \lambda_1^s \varphi_1$ continuous and positive on  $\An$,
for any non-zero
%measure
$\Theta_n$,
\begin{align*}
\alpha_{1,s}
% = \langle \psi_s,\varphi_1\rangle_{L^2(\An)}
=\int_{\An} (K_s \varphi_1)(\ux) \Theta_n(\d \ux)
= \lambda_1^s \int_{\An} \varphi_1 (\ux) \Theta_n(\d \ux)
\end{align*}
is bounded away from zero, uniformly over $s \le \ell+1$. Consequently,
\[ \lim_{T\to\infty} \, 
\frac{\Xi_{n,T}^{(2)}}{\alpha_{1,\cs}^2}  =  \lambda_1  \,.
\]
To treat $\Xi_{n,T}^{(1)}$ of~\eqref{eq:mu-free-3}, note that
by Cauchy--Schwarz and having $\sup_{i \ge 2}
|\widehat \lambda_i| \le 1-\delta$,
\begin{align}\label{eq:err-Xi1}
\bigg|\sum_{\substack{i,j\geq 1 \\ i+j > 2}} \gamma_{i,j}\alpha_{i,\cs} \alpha_{j,\cs}
\widehat\lambda_i^{\ct}\widehat\lambda_j^{\ct-1} \bigg|
&\leq (1-\delta)^{\ct-1} \sum_{i,j}
\left|\gamma_{i,j}\alpha_{i,\cs} \alpha_{j,\cs} \right| \nonumber \\
& \leq (1-\delta)^{\ct-1} \sqrt{\sum_{i,j} \gamma_{i,j}^2}
\sqrt{\sum_{i,j} \alpha_{i,\cs}^2 \alpha_{j,\cs}^2} \nonumber \\
&= (1-\delta)^{\ct-1} \| K^\Gamma_1 \|_{\Ltt} \, \|\psi_{\cs}\|_2^2 \,.
\end{align}
Taking $T\to\infty$, we see that
\[
\lim_{T\to\infty}  \frac{\Xi_{n,T}^{(1)}}{\alpha_{1,\cs}^2} = \gamma_{1,1} \,.
\]
 Altogether, we have thus established that
\begin{equation}
    \label{eq:mu-free-4}
    \lim_{T\to\infty}\mu_{n,T}^\bcfree(\Gamma) = \frac{\gamma_{1,1}}{\lambda_1}\,.
\end{equation}

We now repeat the same analysis for $\mu_{n,T}^\bczero$, where
since the limit as $T \to \infty$ exists, we assume hereafter that
$T$ is integer (and set $s=\ell=1$). Further, for
simplicity we opt to take $\uy=\ux$ and let
$\psi^{(\epsilon)} (\uu) := K_1(\epsilon \ux,\uu)$.
Inferring that $\psi^{(\epsilon)} \in L^2(\An)$ (because
$K_2(\epsilon \ux,\epsilon \ux) < \infty$), we can write
\[
\psi^{(\epsilon)} = \sum_{i=1}^\infty \alpha_i^{(\epsilon)} \varphi_i \,,
\]
where
\[\alpha_i^{(\epsilon)} := \big<\psi^{(\epsilon)},\varphi_i\big>_{L^2(\An)}\,, \quad
\|\psi^{(\epsilon)}\|_{L^2(\An)}^2 = K_2(\epsilon \ux,\epsilon \ux) =
\sum_{i=1}^\infty (\alpha^{(\epsilon)}_i)^2<\infty \,.
\]
The exact same argument then shows that
$\mu_{n,T}^\bczero(\Gamma)$ is the limit at $\epsilon \to 0$ of
$\Xi^{(1,\epsilon)}_{n,T}/\Xi^{(2,\epsilon)}_{n,T}$ where
\begin{equation}
\label{eq:mu-zero-3}
\Xi^{(1,\epsilon)}_{n,T} := \sum_{i,j \ge 1} \gamma_{i,j} \alpha^{(\epsilon)}_i
\alpha^{(\epsilon)}_j \widehat\lambda_i^{T-1} \widehat\lambda_j^{T-2} \,,
\qquad
\Xi^{(2,\epsilon)}_{n,T} :=
\lambda_1\sum_{i=1}^\infty \widehat\lambda_i^{2T-2} (\alpha^{(\epsilon)}_i)^2 \,.
\end{equation}
With $\psi^{(\epsilon)}>0$ and $\varphi_1>0$, we have as before that
$\alpha^{(\epsilon)}_1 > 0$. Moreover, setting
\[ \kappa_\epsilon := \frac{\|\psi^{(\epsilon)}\|^2_{L^2(\An)}}{(\alpha^{(\epsilon)}_1)^2}\,,\]
we have
analogously to~\eqref{eq:err-Xi2} and~\eqref{eq:err-Xi1} that
\begin{align*}
0 \le \frac{\Xi_{n,T}^{(2,\epsilon)}}{(\alpha^{(\epsilon)}_1)^2}
- \lambda_1 & \le \lambda_1 (1-\delta)^{2T-2}
\kappa_\epsilon \,,
\\
\Big|\frac{\Xi_{n,T}^{(1,\epsilon)}}{(\alpha^{(\epsilon)}_1)^2} - \gamma_{1,1} \Big|
& \le  (1-\delta)^{T-2}
\| K^\Gamma_1 \|_{\Ltt} \, \kappa_\epsilon \,.
\end{align*}
We shall employ the following asymptotic as
$\epsilon \to 0$, the proof of which we defer to Section \ref{subsec-2.3}.
\begin{lemma}\label{lem:ker-asymptot}
Setting $\un:=(2n-1,2n-3,\ldots,1)$, we have that
\begin{equation}\label{eq:bd-kap-eps}
\limsup_{\epsilon \to 0} \frac{K_2(\epsilon \un,\epsilon \un)}{
\Big( \int_{u_1 \le 1} K_1(\epsilon \un, \uu) \varphi_1(\uu) d\uu \Big)^2} < \infty \,.
\end{equation}
\end{lemma}
Since $K_1$ and $\varphi_1$ are both positive, \eqref{eq:bd-kap-eps} applies also
without the restriction to $u_1 \le 1$, with Lemma \ref{lem:ker-asymptot} yielding that
$\kappa_\epsilon$ is uniformly bounded (as $\epsilon \to 0$), when $\ux=\un$.  Hence, thanks
to our freedom to choose the boundary, we have that
\begin{equation}
    \label{eq:mu-zero-2}
    \lim_{T\to\infty}\mu_{n,T}^\bczero(\Gamma) = \frac{\gamma_{1,1}}{\lambda_1}\,,
\end{equation}
which in light of~\eqref{eq:mu-free-4} concludes our proof.
\qed

\subsection{Proof of Lemma~\ref{lem:K1-compact}}
Letting
\begin{align*}
    \fB_0 &= \big\{ f \,:\; \|f\|_{L^2(\An)}\leq 1\big\}\qquad\mbox{and} \\
    \fB_1 &= \big\{ (K_1 f)\,:\; f\in \fB_0\big\}\,,
\end{align*}
we will establish compactness by verifying the Fr\'echet--Kolmogorov criteria
(see~\cite[p.~275]{Yosida}, as well as~\cite{Sudakov}).

First, with $\P$ denoting the law of Brownian motion $\{W(t)\}_{t \in [0,1]}$ in $\R^n$ started at the origin and $\E$ its corresponding expectation, note that
\begin{equation}\label{eq:K1f-alt}
(K_1 f)(\ux) = \E \left[ \one_{\Omega_n^{[0,1]}}(\ux + W(\cdot))\, e^{-\cA_{[0,1]}(\ux+W(\cdot))} f(\ux+W(1))\right] \,.
\end{equation}
Now, setting for $f$ supported on $\An$,
\begin{align}\label{eq:M(f)-def}
M(f) & :=  \sup_{\ux\in\An} \E [|f(\ux+W(1))|]
%= \sup_{\ux\in\An} \int_{\An} (2\pi)^{-n/2} e^{-\frac12 \|\uy-\ux\|_2^2} \left|f(\uy)\right|\,\d\uy
\,,
\end{align}
note that by Cauchy--Schwarz,
\begin{equation}\label{eq:M(f)-bound}
M(f)^2 \leq  \sup_{\ux\in\An} \E [f(\ux+W(1))^2] \le
\|f \|^2_{L^2(\An)} \,\sup_{\ux,\uy \in\An} \{\phi_1(\uy-\ux)\}
% \le \|f \|_{L^2(\An)}^2
\leq 1\,,
\end{equation}
where $\phi_v(\cdot)$ denotes the
%Gaussian
density in~\eqref{def:Gauss} and
the last inequality holds for all $f\in\fB_0$.

This readily implies the following uniform bound on $g=K_1f \in \fB_1$,
where by a computation similar to the third line of \eqref{eq-K1K2},
for any $\ux
\in \An$,
\begin{align}\label{eq:g(x)^2-unif-bnd}
    \left|g(\ux)\right| &\leq \int_{\An} \bE_n^{\ux,\uy,[0,1]}\big[\,
    e^{-\fa \int_0^1 X_1(s) ds} \big]|f(\uy)|\,\d\uy
    \leq c e^{-\frac{\fa}{2} x_1}  M(f) \leq
    c e^{-\frac{\fa}{2} x_1}\,,
    \end{align}
for some finite $c=c(\fa)$,
% = \E [exp(-\fa \int_0^1 B_s ds)]
independent of $\ux$ and $f \in \fB_0$.
We deduce in particular that
\begin{equation}\label{eq:equi-tight}
\limsup_{R\to\infty} \sup_{g\in \fB_1} \int_{\substack{\ux\in\An \\ x_1>R}}\left|g(\ux)\right|^2 d\ux =0\,,
\end{equation}
establishing equitightness (and, due to \eqref{eq:g(x)^2-unif-bnd}, also
uniform boundedness, although it is not needed in view of~\cite{Sudakov}).

%Having established equitightness for $\fB_1$,
It remains to establish equicontinuity for $\fB_1$, where in view of
\eqref{eq:equi-tight} and the compactness of $\bAn \cap \{x_1 \le R \}$
it suffices to bound, in terms of $\|\uh\|$,
% and $R$ for $\uh \in \R^n$,
the value of
 \[
 \sup_{g \in \fB_1, \ux \in \An, \ux+\uh \in \An}
\{ |g(\ux+\uh) - g(\ux)|
%\one_{x_1 \le R}
\} \,.
\]
Using the representation \eqref{eq:K1f-alt} for $g=K_1 f$,
%coupling $\bE_n^{\ux+\uh,\uy,[0,1]}$ \blue{and} $\bE_n^{\ux,\uy,[0,1]}$ via
% $\{W(t)\}_{t \in [0,1]}$ in $\R^n$ started at the origin.
%\blue{Specifically, setting}
%\[ g(\ux) := \E \left[ \one_{\Omega_n^{[0,1]}}(\ux + W(\cdot))\,
% e^{-\cA_{[0,1]}(\ux+W(\cdot))} f(\ux+W(1))\right]\,, \]
%we aim to bound $ \| g(\ux+\uh)-g(\ux)\|_{L^2(\An)}$ in terms of $\|\uh\|$.}
%Prior to addressing this quantity, a useful preliminary observation regarding $g(\ux)$ is that the absolute value of the term within the expectation is always at most $|f(\ux+W(1))|$; thus, recalling~\eqref{eq:M(f)-def} and~\eqref{eq:M(f)-bound}, it follows that for every $f\in\fB_0$,
%\begin{equation}\label{eq:g(x)^2-unif-bnd} \sup_{\ux\in\An} g(\ux)^2 \leq \sup_{\ux\in\An} \E[|f(\ux+W(1))|]^2 =
%M(f)^2 \leq \blue{1}\,.
%\end{equation}
we start by reducing to $\tilde g(\cdot)$ in which we extracted out the explicit
dependence of the area tilt on $\ux$. Specifically, let
\[ \tilde g(\ux) := \E \left[ \one_{\Omega_n^{[0,1]}}(\ux + W(\cdot))\, e^{-\cA_{[0,1]}(W(\cdot))} f(\ux+W(1))\right]\,. \]
By a slight abuse of notation, letting $\cA_{[0,1]}(\ux)$ denote $\cA_{[0,1]}(X(\cdot))$ for $X\equiv \ux$, which is nothing but $\fa\left< \underline\fb,\ux\right>$ for $\underline\fb:=(1,\fb,\ldots,\fb^{n-1})$, we see that
\[  g(\ux) = e^{-\cA_{[0,1]}(\ux)} \, \tilde g(\ux)\,,
\]
and therefore,
\begin{align*} | g(\ux)-g(\ux+\uh)|
 &= | e^{-\cA_{[0,1]}(\ux)}\big(\tilde g(\ux)- e^{-\cA_{[0,1]}(\uh)}\tilde g(\ux+\uh)\big)| \\
 &\leq \big| e^{\cA_{[0,1]}(\uh)} - 1\big|\, | g(\ux+\uh)| +
 e^{-\cA_{[0,1]}(\ux)} | \tilde g(\ux)- \tilde g(\ux+\uh) | \,.
 \end{align*}
For the first term note that
$|\cA_{[0,1]}(\uh)| = |\fa \left<\underline\fb,\uh\right>| \leq \fa \|\underline{\fb}\| \|\uh\|$ and though $\uh$ may be outside $\An$,
% $\ux+\uh \in \An$, yielding
by Taylor expansion and~\eqref{eq:g(x)^2-unif-bnd} we have that for any $\|\uh\| \le 1$,
\begin{align*}
\sup_{g\in\fB_1,\ux+\uh \in\An} \big| e^{\cA_{[0,1]}(\uh)}-1\big| \,
| g(\ux+\uh)| \leq C(\fa,\fb,n) \|\uh\|\,.
\end{align*}
Further, with $\cA_{[0,1]}(\ux)\geq 0$ for all $\ux\in\An$,
it remains only to bound $|\tilde g(\ux)-\tilde g(\uy)|$
 uniformly over $g \in \fB_1$, $\ux\in\An$ and $\uy\in \An$
 such that $\|\uy-\ux\| \le \delta$.
To this end, let
\[
\tau_{\ux} := \inf\left\{ t \geq 0 \,:\; \ux+W(t) \notin \An \right\}\,,\quad\mbox{so that }\quad \one_{\Omega_n^{[0,1]}}(\ux+W(\cdot)) = \one_{\left\{ \tau_{\ux} > 1\right\}}\,.
\]
We then have in terms of
 \[ \Delta(\ux,\uy) := \one_{\left\{\tau_{\uy}> 1\right\}} f(W(1)+\uy) - \one_{\left\{\tau_{\ux}> 1\right\}} f(W(1)+\ux)
\]
\and $\eta \in (0,1)$, that
\begin{align*}
| \tilde g(\uy)-\tilde g(\ux)| &=
| \E \big[ e^{-\cA_{[0,1]}(W(\cdot))} \Delta(\ux,\uy) \big]
| \le \E\left[ |\Psi_1|\right] +  | \E\left[\Psi_2\right]|\,,
\end{align*}
where
\begin{align*}
\Psi_1 &:= e^{-\cA^*_{[0,1-\eta]}(W(\cdot))}
\left(
e^{-\cA_{[1-\eta,1]}(W(\cdot)-W(1-\eta))}-1
\right)
\Delta(\ux,\uy)\,, \\
\Psi_2 &:= e^{-\cA^*_{[0,1-\eta]}(W(\cdot))}\Delta(\ux,\uy)\,,
\end{align*}
and
\[ \cA_{[0,1-\eta]}^*(W(\cdot)) = \cA_{[0,1-\eta]}(W(\cdot)) + \cA_{[1-\eta,1]} (W(1-\eta)) \,.
\]
To bound $\E|\Psi_1|$, use the fact that $|\Delta(\ux,\uy)| \leq |f(W(1)+\uy)| + |f(W(1)+\ux)|$ together with H\"older's inequality to infer that
$\E |\Psi_1|$ is at most
\begin{align*}
 \E\bigg[ e^{-4\cA^*_{[0,1-\eta]}(W(\cdot))} \bigg]^{\frac14}
\E\bigg[\left|e^{-\cA_{[1-\eta,1]}(W(\cdot)-W(1-\eta))}-1\right|^4\bigg]^{\frac14} \bigg(2\sup_{\ux\in \An}\E\left[f(W(1)+\ux)^2\right]\bigg)^{\frac12}.
\end{align*}
Noting that the variance of the centered Gaussian $\cA^*_{[0,1-\eta]}(W(\cdot))$ is at most some $v=v(\fa,\fb,n)$ finite,
the first expectation above is uniformly bounded (namely, by $e^{8v}$).
%by Gaussian bounds on $\E \exp(-c\int_0^{1-\eta} B(t) dt)$ for $1$-dimensional Brownian motion $B(t)$; e.g., one can use the fact $|\cA^*_{[0,1-\eta]}(W(\cdot))| \leq \fa \sum_i \fb^{i-1} \sup_{s\leq 1-\eta} |W_i(s)|$.
Similarly, by \eqref{eq:M(f)-bound}, the third term is at most $\sqrt{2}$,
uniformly over $f \in \fB_0$. Finally, with
$\cA_{[1-\eta,1]}(W(\cdot)-W(1-\eta))$
a centered Gaussian of variance $c(\fa,\fb,n) \eta^2$ for some
finite $c(\fa,\fb,n)$,
% \frac{1}{3} \sum_i (\fa \fb^{i-1})^2
%\leq \eta \fa\sum_i \fb^{i-1} \sup_{1-\eta\leq s \leq 1} |W_i(s)-W_i(1-\eta)|$,
the expectation in the second term is at most $\epsilon_0(\eta) \to 0$
as $\eta\to 0$. Overall, we conclude that
\begin{equation}\label{eq:eps1}
\E |\Psi_1| \leq \epsilon_{1}(\eta)\downarrow 0 \quad \hbox{as} \quad
\eta\downarrow 0\,, \quad \hbox{uniformly over} \quad g \in \fB_1, \ux \in \An\,.
\end{equation}
Turning to $\Psi_2=e^{-\cA^*_{[0,1-\eta]}(W(\cdot))}\Delta(\ux,\uy)$, utilizing
 the identity
\[
\one_{\{\tau_{\uy}>1\}} = 1 - \one_{\{\tau_{\uy}\leq 1-\eta \;,\; \tau_{\ux}\leq 1-\eta\}}
-\one_{\{1-\eta  <\tau_{\uy} \leq 1 \}} - \one_{\{\tau_{\uy} \le 1-\eta\;,\; 1-\eta < \tau_{\ux}\leq 1\}} - \one_{\{\tau_{\uy}\leq 1-\eta\:,\: \tau_{\ux}>1\}} \,,
\]
and its dual where the roles of $\tau_{\uy}$ and $\tau_{\ux}$ have been exchanged,
yields the decomposition
\begin{align*}
    \Delta(\ux,\uy) &= \Upsilon_1 - \Upsilon_2(\uy) -\Upsilon_3(\uy,\ux)-\Upsilon_4(\uy,\ux)
    + \Upsilon_2(\ux) + \Upsilon_3(\ux,\uy)+ \Upsilon_4(\ux,\uy) \,,
\end{align*}
where
\begin{align*}
\Upsilon_1 &:= [f(W(1)+\uy) - f(W(1)+\ux)](1-\one_{\{\tau_{\ux}\leq 1-\eta \;,\; \tau_{\uy}\leq 1-\eta\}})\,,\\
\Upsilon_2(\uy) &:= f(W(1)+\uy)\one_{\{1-\eta < \tau_{\uy}\leq 1 \}} \,, \\
    \Upsilon_3(\uy,\ux) &:= f(W(1)+\uy)\one_{\{\tau_{\uy} \le 1-\eta\;,\; 1-\eta < \tau_{\ux}\leq 1\}} \,, \\
    \Upsilon_4(\uy,\ux) &:= f(W(1)+\uy) \one_{\{\tau_{\uy}\leq 1-\eta\:,\: \tau_{\ux}>1\}} \,.
\end{align*}
For the contribution to $|\E[\Psi_2]|$
%from the term that corresponds
due to $\Upsilon_1$, condition on
$\cF_{1-\eta} = \sigma(\{W(s)\}_{s\leq 1-\eta})$,
on which the indicator in $\Upsilon_1$ is measurable, to get
\begin{align*}
     \left|\E\big[ \Upsilon_1 e^{-\cA^*_{[0,1-\eta]}(W(\cdot))}\big]\right| &\leq \E\left[e^{-\cA^*_{[0,1-\eta]}W(\cdot))}\right] \\
     &\cdot \sup_{\uz}
\Big| \E\big[f(W(1)+\uy)- f(W(1)+\ux) \;\big|\; W(1-\eta)=\uz\big]\Big|\,.
\end{align*}
While treating $\E|\Psi_1|$, we saw that the first term on the right-hand
is some finite $C(\fa,\fb,n)$, independently of $\ux,\uh$.
For the second term, extending $f\in\fB_0$ from $\An$ to $\R^n$
via $f(\ux)=0$ for $\ux\notin \An$, yields that
$\|f\|_2 = \|f\|_{L^2(\An)} \leq 1$. Thus, performing
a change of variable
 $\uv := W(1)+\uy$ in $\E[f(W(1)+\uy) \mid W(1-\eta)=\uz]$ and
 $\uv := W(1)+\ux$ in $ \E[f(W(1)+\ux)\mid W(1-\eta)=\uz]$, we
 get that the absolute difference between these expectations is
\begin{align*}
   & \Big| \int \big[\phi_\eta(\uv-\uy-\uz) -\phi_\eta(\uv-\ux-\uz)\big] f(\uv)\,\d\uv \Big| \\
    & \leq \|f\|_{L^2(\An)} \eta^{-n/4}
    \big\|  \phi_1(\uw-\eta^{-1/2} (\uy-\ux))-\phi_1(\uw)
    \big\|_2
    \le C(n) \, \eta^{-n/4-1/2} \delta   \,,
\end{align*}
where the first inequality is obtained by Cauchy--Schwarz and an
additional change of variable $\uw = \eta^{-1/2} (\uv-\uz-\ux)$, and
the second inequality by an easy computation (utilizing that $1-e^{-r} \le r$).
% with C(n)^2 = 2 \pi^{n/2} (2\pi)^{-n} \frac{1}{4} \le 1.
%
Thus, choosing
\begin{equation}\label{eq:eta-def}
\delta \le \eta^{n/4+1},
\end{equation}
makes the contribution of $\Upsilon_1$ negligible.
% (noting that we have flexibility in the choice of $\eta$ as a function of $\|\uh\|$)

To deal with the contribution of $\Upsilon_2(\uy)$ to $|\E[\Psi_2]|$, observe that
by H\"older's inequality,
\begin{align}\label{eq:basic-holder}
\E&\left[
e^{-\cA^*_{[0,1-\eta]}(W(\cdot))}
\one_{\{1-\eta < \tau_{\uy} \leq 1\}}
|f(W(1)+\uy)|
\right] \nonumber \\
& \leq
 \E\Big[ e^{-4\cA^*_{[0,1-\eta]}(W(\cdot))} \Big]^{\frac14}
\P\Big(1-\eta < \tau_{\uy} \leq 1\Big)^{\frac14}
\Big(\sup_{\uy\in \An}\E\left[f(W(1)+\uy)^2\right]\Big)^{\frac12}.
\end{align}
While bounding $\E|\Psi_1|$ we have seen that the first and
third terms are at most some $c(\fa,\fb,n)$ finite, uniformly over $\fB_0$, so
it suffices to show that
\begin{equation}\label{eq:eps2}
\epsilon_2(\eta):=\sup_{\uy \in \An} \{
\P(1-\eta < \tau_{\uy}\leq 1) \} \to 0 \qquad \mbox{as \;\; $\eta\to 0$}\,.
\end{equation}
Indeed, taking a union bound over the $n$ different boundaries of $\An$ that are considered in $\tau_{\uy}$, reduces, up to the factor $n$, to the bound in case $n=1$, namely for the first hitting time $T_b$
of level $-b<0$ by a standard Brownian motion $B_t$. The corresponding probability density
$f_{T_b}(t)=be^{-b^2/(2t)}/\sqrt{2\pi t^3}$
% the latter hitting time
is bounded, uniformly over $b$ and $t \ge 1/2$,
% See (10.1.4) in 310 lecture notes.
thereby yielding \eqref{eq:eps2}.

The same analysis applies to the contributions
%from $\Upsilon_2(\ux)$ and
from the $\Upsilon_3$ terms.

Analogously to \eqref{eq:basic-holder} the contribution of
$\Upsilon_4(\uy,\ux)$ to $|\E[\Psi_2]|$ is bounded above by
\begin{align*}
 & \E\left[e^{-\cA_{[0,1-\eta]}^*(W(\cdot))}
  \one_{\{\tau_{\uy}\leq 1-\eta \;,\; \tau_{\ux}> 1\}} |f(W(1)+\uy)|
  \right]  \\
 \quad&\leq C \, \|f\|_{L^2(\An)}
 \P(\tau_{\uy}\leq 1-\eta\;, \;\tau_{\ux} > 1)^{1/4} \le
 C \, \epsilon_3(\delta,\eta)^{1/4} \,,
 \end{align*}
for some $C=C(\fa,\fb,n)$, any $f \in \fB_0$ and
\[
\epsilon_3(\delta,\eta) := \sup_{\ux,\uy \in \An, \|\ux-\uy\| \le \delta}
\P(\tau_{\uy}\leq 1-\eta\;, \;\tau_{\ux} > 1) \,.
\]
With the same bound applying for $\Upsilon_4(\ux,\uy)$, it remains only
to show that $\epsilon_3(\delta,\eta) \to 0$ as $\delta \to 0$ (for any
fixed $\eta>0$). To this end, by a union bound over the $n$ different boundaries of $\An$, as done for proving \eqref{eq:eps2},
the probability in question is at most~$n$ times the probability
that standard Brownian motion $B(t):=\frac{1}{\sqrt{2}}(W_i(t)-W_{i+1}(t))$
reach level $-b$ by time $1-\eta$ (here $b=(y_i-y_{i+1})/\sqrt{2}$),
while remaining above $-(b+\delta)$
% since (x_i-x_{i+1})-(y_i-y_{i+1}) \le \delta \sqrt{2}.
up till time $1$. With Brownian motion a strong Markov process
of independent increments, we thus deduce by the reflection principle
that
\begin{align*}
n^{-1} \epsilon_3(\delta,\eta) \le \P(\inf_{s\leq \eta} \{B(s)\} > -\delta)
 = 1 - 2\P(B(\eta) \ge \delta) = \P(|B(\eta)| <\delta)\,,
 \end{align*}
which goes to zero as $\delta \to 0$ (for any fixed $\eta>0$).
\qed

\subsection{Proof of Lemma \ref{lem:ker-asymptot}}\label{subsec-2.3}

Setting $\widehat K_t$ for the operator $K_t$ in the case $\fa=0$ (no area tilt),
we first establish~\eqref{eq:bd-kap-eps} for $\widehat K_t$. Namely, we show
that,
\begin{equation}\label{eq:bd-kap-a}
\limsup_{\epsilon \to 0} \frac{\widehat K_2(\epsilon \un,\epsilon \un)}{
\Big( \int_{u_1 \le 1} \widehat K_1(\epsilon \un, \uu) \varphi_1(\uu) d\uu \Big)^2} < \infty \,.
\end{equation}
Our starting point for \eqref{eq:bd-kap-a} is the following explicit formula,
valid for any $\uy \in \An$ and any $t,\epsilon>0$,
\begin{equation}\label{eq:exact}
\widehat K_t(\epsilon \un,\uy) = 2^{n^2} \phi_t(\uy)
e^{-\epsilon^2\|\un\|^2/(2t)} \prod_{i} \sinh\big(\frac{\epsilon y_i}{t}\big)
\prod_{j<k} \Big[
\sinh^2\big(\frac{\epsilon y_j}{t}\big) - \sinh^2\big(\frac{\epsilon y_k}{t}\big) \Big] \,.
\end{equation}
Indeed, for $\epsilon=1$ this is the explicit evaluation in
\cite[Display below (24)]{Grabiner99} of the Karlin--McGregor
determinantal formula~\cite{KarlinMcGregor59} for the transition kernel,
\[
q_t(x,y)=\phi_t(y-x)-\phi_t(y+x)=2 \phi_t(y) e^{-x^2/(2t)} \sinh(xy/t) \,,
\]
of a scalar Brownian motion absorbed at level zero, when starting at the
distinguished point $\un$. We thus get \eqref{eq:exact} by noting
that the non-trivial factors $\sinh(x_i y_j/t)$ are invariant to
changing from $(\epsilon \un,\uy)$ to $(\un,\epsilon \uy)$.

In particular, with $g(x):=\sinh(x/2)$ being zero at $x=0$
and globally Lipschitz($L$) on $[0,2n]$, we get from \eqref{eq:exact} that
for some $c_n$, $C_n$ finite and any $\epsilon \in [0,1]$,
\begin{align}\label{eq:ubd-wK2}
\widehat K_2(\epsilon \un,\epsilon \un) &\le c_n
\prod_i g(\epsilon^2 n_i) \prod_{j<k} [g^2(\epsilon^2 n_j) - g^2(\epsilon^2 n_k)]
\nonumber \\
&\le c_n L^{n^2} \prod_i (\epsilon^2 n_i) \prod_{j<k} \big[(\epsilon^2 n_j)^2
- (\epsilon^2 n_k)^2 \big] = C_n \epsilon^{2 n^2} \,.
\end{align}
Next, noting that on $\R_+$ both $\sinh(x) \ge x$ and $\sinh^2(x)-x^2$ are
non-decreasing, we deduce from \eqref{eq:exact} that for any $\uu \in \An$
and $\epsilon \in [0,1]$,
\[
\widehat K_1(\epsilon \un,\uu) \ge 2^{n^2} e^{-\|\un\|^2/2} \epsilon^{n^2}
\hat{\phi}(\uu) \,, \quad \hbox{where} \quad
\hat{\phi}(\uu) := \phi_1(\uu) \prod_{i} u_i \prod_{j<k} (u_j^2 - u_k^2) \,.
\]
With $\hat{\phi}(\cdot)$ and $\varphi_1(\cdot)$ positive
on $\An$, we get from the latter bound that
\[
\inf_{\epsilon \in [0,1]} \epsilon^{-n^2}
\int_{u_1 \le 1} \widehat K_1(\epsilon \un, \uu) \varphi_1(\uu) d\uu > 0 \,,
\]
which in combination with \eqref{eq:ubd-wK2} establishes~\eqref{eq:bd-kap-a}.
% To obtain its analog for $K_t$, i.e., the sought bound~\eqref{eq:bd-kap-eps},

Next, recall that $K_t(\ux,\uy)$ is point-wise decreasing in $\fa$ and
in particular bounded from above by $\widehat K_t(\ux,\uy)$;
%since the tilt term $e^{-\cA_{I}(X(\cdot))}$ in $K_t$ (missing from $\widehat K_t$) is at most $1$ for every $X\in\Omega_n^I$;
thus, the sought bound \eqref{eq:bd-kap-eps} for $K_t$ follows from~\eqref{eq:bd-kap-a} once we show that for some finite $C=C(\fa,\fb,n)$
and any $\uu \in \An$ with $u_1\leq 1$,
\begin{equation}\label{eq:unif-bd}
\sup_{\epsilon\in (0,1]} \Big\{ \frac{\widehat K_1(\epsilon\un,\uu)}{K_1(\epsilon\un,\uu)} \Big\} \le C \,.
\end{equation}
Turning to the latter bound, we define for finite $M$ the event
\[ \Gamma_M := \left\{ \max_{t \in [0,1]} \{X_1(t)\} \le M \right\}\,,\]
noting that for $c := \fa \langle \underline{\fb}, \underline 1 \rangle$,
any $u_1 \le 1$ and $\epsilon \le 1$,
\begin{align*}
K_1(\epsilon\un,\uu) &\geq e^{-c M} \bE_n^{\epsilon\un,\uu,[0,1]}
\left[ \one_{\Gamma_M} \one_{\Omega_n^{[0,1]}}\right]
%= e^{-M} \widehat K_1^{\Gamma_M}(\epsilon\un,\epsilon\uu)\\
= e^{-c M} \widehat K_1(\epsilon\un,\uu) \widehat\P^{\epsilon\un,\uu,[0,1]}_n (\Gamma_M) \\
&\geq e^{-c M} \widehat K_1(\epsilon\un,\uu) \widehat\P^{\un,\un,[0,1]}_n (\Gamma_M) \,,
\end{align*}
where $\widehat\P^{\ux,\uy,[0,1]}_n$ is the measure $\P^{\ux,\uy,[0,1]}_n$ from~\eqref{eq:P-ux-uy-def} corresponding to $\fa=0$, and with the second inequality due to~\cite[Lemma~2.7]{CorwinHammond14} (taking there $A=[0,1]$,
$f \equiv 0$,
%); indeed, conditional on
%the event \blue{$\Omega_n^{[0,1]}$} that the
%lines in our ensemble do not intersect each other, nor the line $f\equiv 0$,
%the aforementioned lemma shows that the ensemble of $n$ Brownian bridges on $[0,1]$ with endpoints $\ux,\uy$
%conditioned on $\cE_n$,
%stochastically dominates the analogous ensemble with endpoints $\ux',\uy'$ that satisfy $\ux \geq \ux'$ and $\uy \geq \uy'$ point-wise (
noting that $\un > \uu$ and $\un > \epsilon \un$
whenever $u_1 \leq 1$ and $\epsilon \le 1$
and that the event $\Gamma_M$ is decreasing).

Finally, moving to the unconditional space of $n$ independent bridges rooted at $\un,\un$ via a multiplicative cost of at most $1/\widehat K_1(\un,\un)$, we see that
$\widehat\P^{\un,\un,[0,1]}_n (\Gamma_M^c)$ is at most $\P(\sup_{s \in [0,1]}
\{B(s)\} > M-2n) /\widehat K_1(\un,\un)$
for a one dimensional Brownian bridge from $(0,1)$ to $(1,1)$. By the tightness of the maximum of the latter bridge (and recalling that
$\widehat{K}_1(\un,\un)>0$), one thus has for $M$ large,
depending only on $n$, that
\[
\widehat\P_n^{\un,\un,[0,1]}(\Gamma_M) \ge \tfrac12\,.
\]
Combining the last two displays yields \eqref{eq:unif-bd},
thereby completing the proof.
\qed

\subsection*{Acknowledgment}
We thank Ivan Corwin for bringing to our attention the paper~\cite{Grabiner99}.
We thank the referees for a careful reading of the manuscript and their comments.
A.D.\ was supported in part by NSF grant DMS-1954337.
E.L.\ was supported by NSF grants DMS-1812095 and DMS-2054833.
O.Z.\ was partially supported by the European Research Council (ERC) under the European Union's Horizon 2020 research and innovation programme (grant agreement No.~692452).
This research was further supported in part by BSF grant 2018088.
\bibliographystyle{abbrv}
%\bibliography{bm_tilt}

\begin{thebibliography}{10}

\bibitem{BEF86}
J.~Bricmont, A.~El~Mellouki, and J.~Fr\"{o}hlich.
\newblock Random surfaces in statistical mechanics: roughening, rounding,
  wetting,{$\ldots\,$}.
\newblock {\em J. Statist. Phys.}, 42(5-6):743--798, 1986.

\bibitem{CIW19a}
P.~Caputo, D.~Ioffe, and V.~Wachtel.
\newblock Confinement of {B}rownian polymers under geometric area tilts.
\newblock {\em Electron. J. Probab.}, 24:Paper No. 37, 21, 2019.

\bibitem{CIW19b}
P.~Caputo, D.~Ioffe, and V.~Wachtel.
\newblock Tightness and line ensembles for {B}rownian polymers under geometric
  area tilts.
\newblock In {\em Statistical mechanics of classical and disordered systems},
  volume 293 of {\em Springer Proc. Math. Stat.}, pages 241--266. Springer,
  Cham, 2019.

\bibitem{CLMST12}
P.~Caputo, E.~Lubetzky, F.~Martinelli, A.~Sly, and F.~L. Toninelli.
\newblock The shape of the {$(2+1)D$} {SOS} surface above a wall.
\newblock {\em C. R. Math. Acad. Sci. Paris}, 350(13-14):703--706, 2012.

\bibitem{CLMST14}
P.~Caputo, E.~Lubetzky, F.~Martinelli, A.~Sly, and F.~L. Toninelli.
\newblock Dynamics of {$(2+1)$}-dimensional {SOS} surfaces above a wall: {S}low
  mixing induced by entropic repulsion.
\newblock {\em Ann. Probab.}, 42(4):1516--1589, 2014.

\bibitem{CLMST16}
P.~Caputo, E.~Lubetzky, F.~Martinelli, A.~Sly, and F.~L. Toninelli.
\newblock Scaling limit and cube-root fluctuations in {SOS} surfaces above a
  wall.
\newblock {\em J. Eur. Math. Soc. (JEMS)}, 18(5):931--995, 2016.

\bibitem{CorwinHammond14}
I.~Corwin and A.~Hammond.
\newblock Brownian {G}ibbs property for {A}iry line ensembles.
\newblock {\em Invent. Math.}, 195(2):441--508, 2014.

\bibitem{FS05}
P.~L. Ferrari and H.~Spohn.
\newblock Constrained {B}rownian motion: fluctuations away from circular and
  parabolic barriers.
\newblock {\em Ann. Probab.}, 33(4):1302--1325, 2005.

\bibitem{Fukushima}
M.~Fukushima, Y.~Oshima, and M.~Takeda.
\newblock {\em Dirichlet forms and symmetric Markov processes}.
\newblock De Gruyter, Berlin, second edition, 2011.

\bibitem{Grabiner99}
D.~J. Grabiner.
\newblock Brownian motion in a {W}eyl chamber, non-colliding particles, and
  random matrices.
\newblock {\em Ann. Inst. H. Poincar\'{e} Probab. Statist.}, 35(2):177--204,
  1999.

\bibitem{IOSV21}
D.~Ioffe, S.~Ott, S.~Shlosman, and Y.~Velenik.
\newblock Critical prewetting in the 2{D} {I}sing model.
\newblock {\em Ann. Probab.}, 50(3):1127--1172, 2022.

\bibitem{ISV15}
D.~Ioffe, S.~Shlosman, and Y.~Velenik.
\newblock An invariance principle to {F}errari-{S}pohn diffusions.
\newblock {\em Comm. Math. Phys.}, 336(2):905--932, 2015.

\bibitem{IV18}
D.~Ioffe and Y.~Velenik.
\newblock Low-temperature interfaces: prewetting, layering, faceting and
  {F}errari-{S}pohn diffusions.
\newblock {\em Markov Process. Related Fields}, 24(3):487--537, 2018.

\bibitem{IVW18}
D.~Ioffe, Y.~Velenik, and V.~Wachtel.
\newblock Dyson {F}errari-{S}pohn diffusions and ordered walks under area
  tilts.
\newblock {\em Probab. Theory Related Fields}, 170(1-2):11--47, 2018.

\bibitem{KarlinMcGregor59}
S.~Karlin and J.~McGregor.
\newblock Coincidence probabilities.
\newblock {\em Pacific J. Math.}, 9:1141--1164, 1959.

\bibitem{LMS16}
E.~Lubetzky, F.~Martinelli, and A.~Sly.
\newblock Harmonic pinnacles in the discrete {G}aussian model.
\newblock {\em Comm. Math. Phys.}, 344(3):673--717, 2016.

\bibitem{MZ16}
P.~Maillard and O.~Zeitouni.
\newblock Slowdown in branching {B}rownian motion with inhomogeneous variance.
\newblock {\em Ann. Inst. Henri Poincar\'{e} Probab. Stat.}, 52(3):1144--1160,
  2016.

\bibitem{ReedSimonI}
M.~Reed and B.~Simon.
\newblock {\em Methods of modern mathematical physics. {I}. {F}unctional
  analysis}.
\newblock Academic Press, New York-London, 1972.

\bibitem{Sudakov}
V.~N. Sudakov.
\newblock Criteria of compactness in function spaces.
\newblock {\em Uspekhi Mat. Nauk}, 12:221--224, 1957.

\bibitem{Yosida}
K.~Yosida.
\newblock {\em Functional Analysis}.
\newblock Springer-Verlag, Heidelberg, sixth edition, 1980.

\end{thebibliography}

\end{document}